\documentstyle[amssymb,amsfonts,12pt,psfig]{article}

\newcommand{\bc}{\begin{center}}
\newcommand{\ec}{  \end{center}}

\newcommand{\bi}{\begin{itemize}}
\newcommand{\ei}{  \end{itemize}}
\newcommand{\benum}{\begin{enumerate}}
\newcommand{\eenum}{  \end{enumerate}}
\newcommand{\bd}{\begin{description}}
\newcommand{\ed}{  \end{description}}

\def\R{{\hbox{\bf R}}}

\def\Z{{\hbox{\bf Z}}}

\def\eps{{\varepsilon}}
\def\implies{{\Rightarrow}}

\def\emph#1{{\it #1}}

\begin{document}

\title{From rotating needles to stability of waves: emerging connections between combinatorics, analysis, and PDE.}

\maketitle

\section{Introduction}

In 1917 S. Kakeya posed the {\it Kakeya needle problem}: what is the smallest area which is required to rotate a unit line segment (a ``needle'') by 180 degrees in the plane?  Rotating around the midpoint requires $\pi/4$ units of area, whereas a ``three-point U-turn'' requires $\pi/8$.  In 1927 the problem was answered by A. Besicovitch, who gave the surprising answer that one could rotate a needle using arbitrarily small area.

At first glance, Kakeya's problem and Besicovitch's resolution appear to be little more than mathematical curiosities.  However, in the last three decades it has gradually been realized that this type of problem is connected to many other, seemingly unrelated problems in number theory, geometric combinatorics, arithmetic combinatorics, oscillatory integrals, and even the analysis of dispersive and wave equations.  

The purpose of this article is to discuss the interconnections between these fields, with an emphasis on the connection with oscillatory integrals and PDE.  Two previous surveys (\cite{Wo} and \cite{Bo}) have focused on the connections between Kakeya-type problems with other problems in discrete combinatorics and number theory. 

These areas are very active, but despite much recent progress, our understanding of the problems and their relationships to each other is far from complete.  Ideas from other fields may well be needed to make substantial new breakthroughs. 

\section{Kakeya type problems}

Besicovitch's solution to the Kakeya needle problem relied on two observations.  The first observation, which is elementary, is that one can translate a needle to any location using arbitrarily small area.  The second observation is that one can construct open subsets of $\R^2$ of arbitrarily small area which contain a unit line segment in every direction.  A typical way to construct such sets (not Besicovitch's original construction) is sketched in Figure 2; for a more detailed construction, see \cite{Wo}.

\begin{figure}[htbp] \centering
\ \psfig{figure=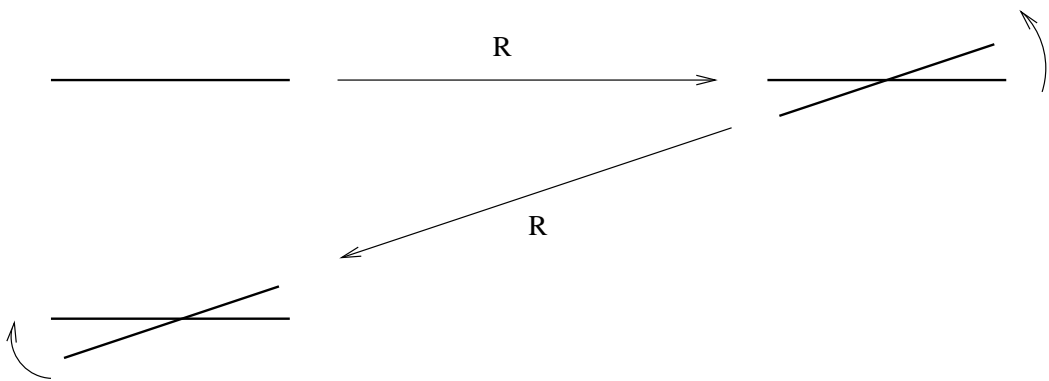}
\caption{To translate a needle, slide it by $R$ units, rotate by roughly $1/R$,
slide it back, and rotate back.  This costs $O(1/R)$ units of area,
where $R$ is arbitrary.
 }
\end{figure}

\begin{figure}[htbp] \centering
\ \psfig{figure=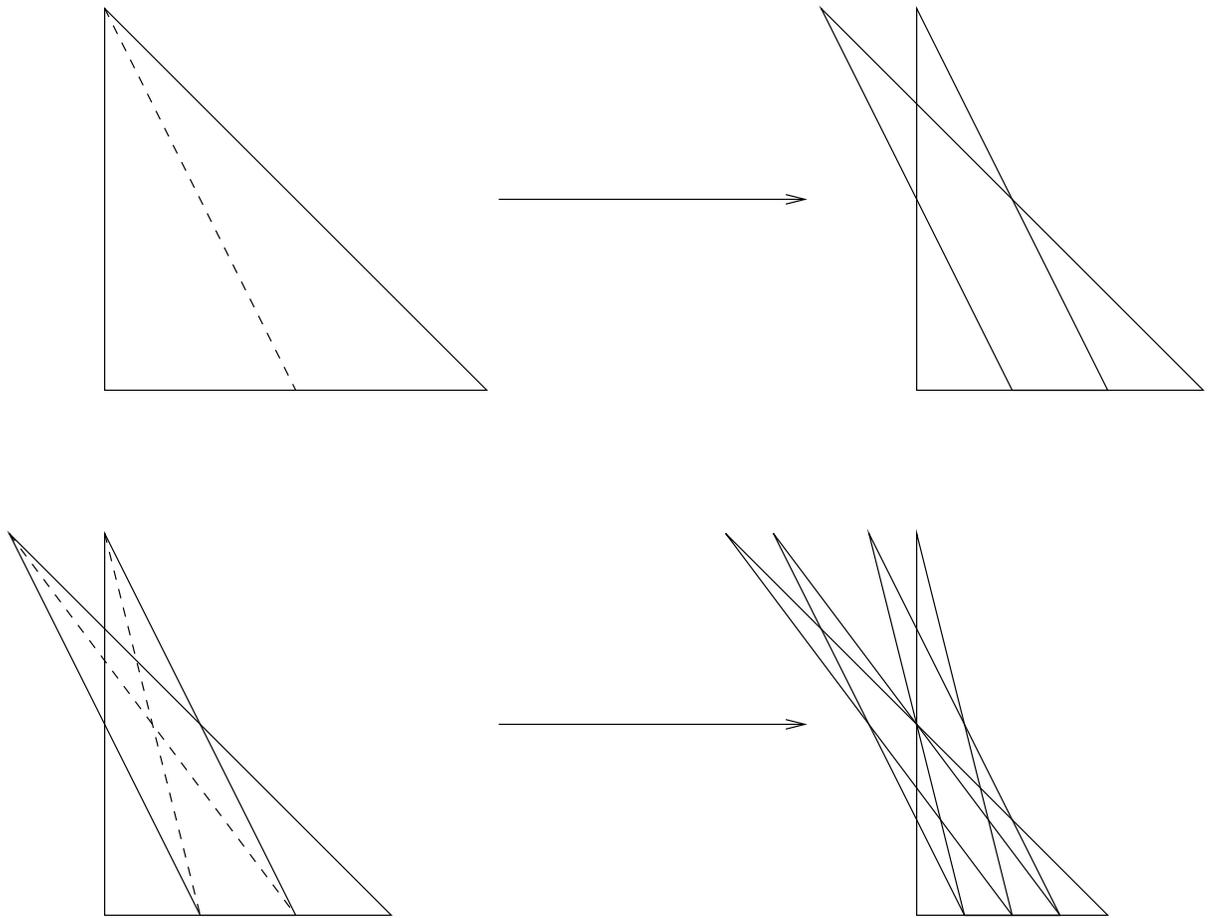}
\caption{The iterative construction of a Besicovitch set.  Each stage consists
of the union of triangles.  To pass to the next stage, the triangles are bisected and shifted together to decrease their area.
 }
\end{figure}

For any $n \geq 2$, define an {\it Besicovitch set} to be a subset of $\R^n$ which contains a unit line segment in every direction.  The construction of Besicovitch shows that such sets can have arbitrarily small measure in any dimension, and can even be made to be measure zero.  Intuitively, this states that it is possible to compress a large number of non-parallel unit line segments into an arbitrarily small set.  

In applications one wishes to obtain more quantitative understanding of this compression effect by introducing a spatial discretization.  For instance, one could replace unit line segments by $1 \times \delta$ tubes for some $0 < \delta \ll 1$ and ask for the optimal compression of these tubes.  Equivalently, one can ask for bounds of the volume of the $\delta$-neighbourhood of a Besicovitch set.

Rather surprisingly, these bounds are logarithmic in two dimensions.  It is known that the $\delta$-neighbourhood of a Besicovitch set in $\R^2$ must have area at least\footnote{Throughout this article, the letter $C$ will denote various constants which vary from line to line.} $C/\log(1/\delta)$; this basically follows from the geometric observation that the area of the intersection of two $1 \times \delta$ rectangles varies inversely with the angle between the long axes of the rectangles.  Recently, U. Keich has shown that this bound is sharp.

This observation can be rephrased in terms of the {\it Minkowski dimension} of the Besicovitch set.  Recall that a bounded set $E$ has Minkowski dimension $\alpha$ or less if and only if for every $0 < \delta \ll 1$ and $0 < \eps \ll 1$, one can cover $E$ by at most $C_\eps \delta^{-\alpha+\eps}$ balls of radius $\delta$.  From the previous discussion we thus see that Besicovitch sets in the plane must have Minkowski dimension 2.

the analogous question in $\R^n$ is known as the {\it Kakeya conjecture}:

\centerline{Does every Besicovitch set in $\R^n$ have Minkowski dimension $n$?}

Equivalently, the Kakeya conjecture asserts that the volume of the $\delta$-neighbourhood of a Besicovitch set in $\R^n$ is bounded below by $C_{n,\eps} \delta^\eps$ for any $\eps > 0$ and $0 < \delta \ll 1$.

There are several variants of this conjecture (e.g. one could discuss Hausdorff dimension instead of Minkowski dimension), but we shall not discuss these variants here for sake of exposition.

The Kakeya conjecture is remarkably difficult.  It remains open in three and higher dimensions, although rapid progress has been made in the last few years.  The best known lower bounds for the Minkowski dimension at this time of writing is
$$ \max(\frac{n+2}{2} + 10^{-10}, \frac{4n+3}{7}),$$
although I expect further improvements to follow very soon.

One can discretize the conjecture.  Let $\Omega$ be a maximal $\delta$-separated subset of the sphere $S^{n-1}$ (so that $\Omega$ has cardinality approximately $\delta^{1-n}$), and for each $\omega \in \Omega$ let $T_\omega$ be a $\delta \times 1$ tube oriented in the direction $\omega$.  The Kakeya conjecture then asserts logarithmic-type lower bounds on the quantity $|\bigcup_{\omega \in \Omega} T_\omega|$.

The above formulation is reminiscent of existing results in combinatorics concerning the number of incidences between lines and points, although a formal connection cannot be made because the nature of the intersection of two $\delta \times 1$ tubes depends on the angle between the tubes, whereas the intersection of two lines is a point regardless of what angle the lines make.  However, it is plausible that one can use the ideas from combinatorial incidence geometry to obtain progress on this problem.  For instance, it is fairly straightforward to show the Minkowski dimension of Besicovitch sets is at least $(n+1)/2$ purely by using the fact that given any two points that are a distance roughly 1 apart, there is essentially only one $\delta \times 1$ tube which can contain them both.

In the 1990s, work by J. Bourgain, T. Wolff, W. Schlag, A. Vargas, N. Katz, I. Laba, the author, and others pushed these ideas further.  For instance, the lower bound of $(n+2)/2$ for the Minkowski dimension was shown in 1995 by Wolff and relies on the $\delta$-discretized version of the geometric statements that every non-degenerate triangle lies in a unique two-dimensional plane, and every such plane contains only a one-parameter set of directions.  However, there appears to be a limit to what can be achieved purely by applying elementary incidence geometry facts and standard combinatorial tools (such as those from extremal graph theory).  More sophisticated geometric analysis seems to reveal that a counterexample to the Kakeya conjecture, if it exists, must have certain rigid structural properties (for instance, the line segments through any given point should all lie in a hyperplane).  These type of ideas have led to a very small recent improvement in the Minkowski bound to $(n+2)/2 + 10^{-10}$, but they are clearly insufficient to resolve the full conjecture.

The Kakeya problem is a representative member of a much larger family of problems of a similar flavour (but with more technical formulations).  For instance, one can define a $\beta$-set to be a subset of the plane which contains a $\beta$-dimensional subset of a unit line segment in every direction.  It is then an open problem to determine, for given $\beta$, the smallest possible dimension of a $\beta$-set.  Low-dimensional examples of such sets arise in the work of H. Furstenburg, and it seems that one needs to understand these generalizations of Besicovitch sets in order to fully exploit the connection between Kakeya problems and oscillatory integrals, which we discuss below.  Other variants include replacing line segments by circles or light rays, considering finite geometry analogues of these problems, or by replacing the quantity $|\bigcup_{\omega} T_\omega|$ by the variant $\| \sum_\omega \chi_{T_\omega} \|_p$ (the relevant conjecture here is known as the {\it Kakeya maximal function conjecture}).  Another interesting member of this family is the {\it Falconer distance set conjecture}, which asserts that whenever $E$ is a compact one-dimensional subset of $\R^2$, that the {\it distance set} $\{ |x-y|: x,y \in E\}$ is a one-dimensional subset of $\R$.  The discrete version of this is the {\it Erd\"os distance problem} - what is the least number of distances determined by $n$ points? - and is also unsolved.  For a thorough survey of most of these questions we refer the reader to \cite{Wo}; see also \cite{kt:falc-furst}.  

\section{The connection with arithmetic combinatorics.}

The Kakeya problem looks like very geometrical, and it is natural to apply elementary incidence geometry facts to bear on this problem.  Although this approach has had some success, it does not seem sufficient to solve the problem.

In 1998 Bourgain introduced a new type of argument, based on {\it arithmetic} combinatorics (the combinatorics of sums and differences), which gave improved results on this problem, especially in high dimensions.  The connection between Kakeya problems and the combinatorics of addition can already be seen by considering the analogy between line segments and arithmetic progressions.  (Indeed, the Kakeya conjecture can be reformulated in terms of arithmetic progressions, and this can be used to connect the Kakeya conjecture to several difficult conjectures in number theory such as the Montgomery conjectures for generic Dirichlet series.  We will not discuss this connection here, but refer the reader to \cite{Bo}).

Bourgain's argument relies on the following ``three-slice'' idea.  Let $\Omega$ and $T_\omega$ be as in the previous section.  We may assume that the tubes $T_\omega$ are contained in a fixed ball.  Suppose that $|\bigcup_{\omega \in \Omega} T_\omega|$ is comparable to $\delta^\alpha$ for some constant $\alpha$; our objective is to give upper bounds on $\alpha$, and eventually to show that $\alpha$ must be zero.

By choosing an appropriate set of co-ordinates, one can ensure that the three slices
$$ X_t := \{ x \in \R^{n-1}: (x,t) \in \bigcup_{\omega in \Omega} T_\omega \}$$
have measure comparable to $\delta^\alpha$ for $t = 0, 1/2, 1$.  Because of the $\delta$-discretized nature of the problem, one can also assume that the discrete set
$$ A_t := X_t \cap \delta \Z^{n-1}$$
has cardinality comparable to $\delta^{\alpha+1-n}$ for $t=0, 1/2, 1$.

\begin{figure}[htbp] \centering
\ \psfig{figure=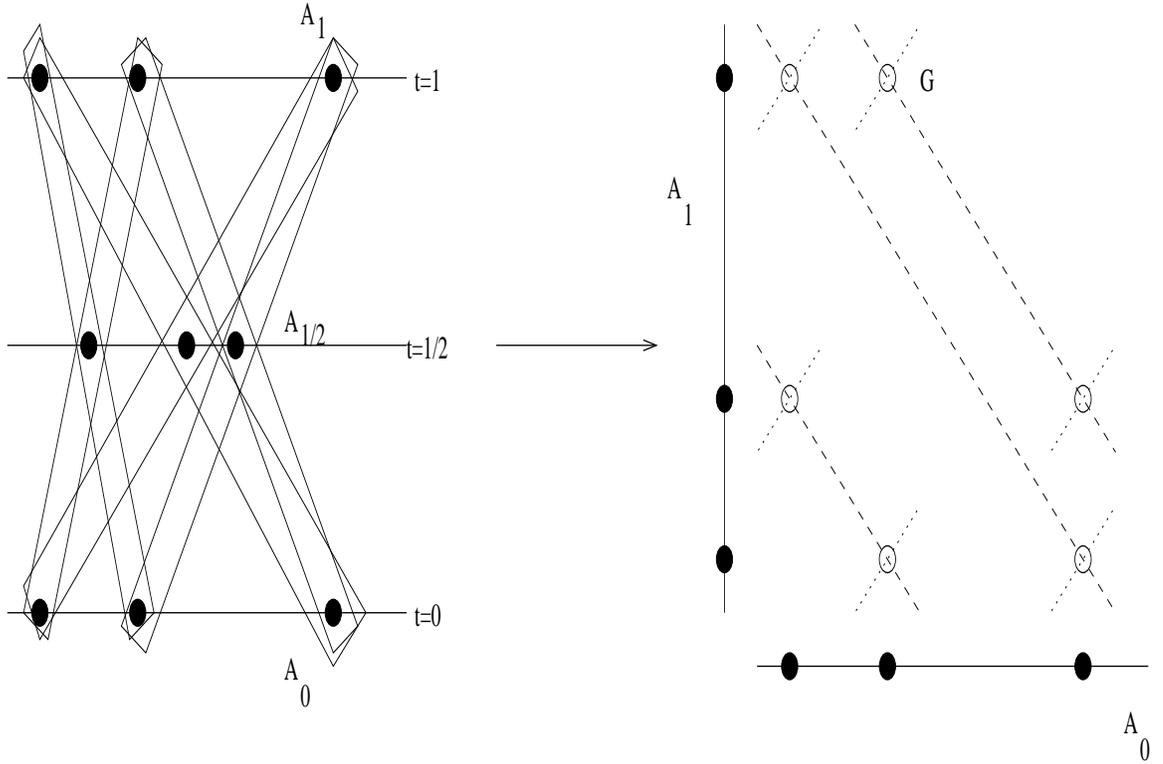,height=4in,width=6in}
\caption{The left picture depicts six tubes $T_\omega$ pointing in 
different directions, and the three
discretized slices $A_0$, $A_{1/2}$, $A_1$, which in this case are three-element
sets.  The right picture depicts the set $G$ associated to this collection of
tubes.  Note that the maps $(a,b) \to a$, $(a,b) \to b$, $(a,b) \to a+b$
have small range (mapping to $A_0$, $A_1$, and $2 A_{1/2}$ respectively),
but the map $(a,b) \to a-b$ is one-to-one.
 }
\end{figure}

Morally speaking, every tube $T_\omega$ intersects each of the sets $A_0$, $A_{1/2}$, $A_1$ in exactly one point.  Assuming this, we see that every tube $T_\omega$ is associated with an element of $A_0 \times A_1$.  Because two points determine a line, these elements are essentially disjoint as $\omega$ varies.  Let $G$ denote the set of all pairs of $A_0 \times A_1$ obtained this way.  Thus $G$ has cardinality about $\delta^{1-n}$.

The sum-set
$$ \{ a+b: (a,b) \in G\}$$
of $G$ is essentially contained inside a dilate of the set $A_{1/2}$; this reflects the fact that the intersection of $T_\omega$ with $A_{1/2}$ is essentially the midpoint of the intersection of $T_\omega$ with $A_0$ and $A_1$.  In particular, the sum-set is quite small, having cardinality only $\delta^{\alpha+1-n}$.  On the other hand, the difference set
$$ \{a-b: (a,b) \in G\}$$
of $G$ is quite large, because the tubes $T_\omega$ all point in different directions.  Indeed, this set has the same cardinality as $G$, i.e. about $\delta^{1-n}$.

Thus, if $\alpha$ is non-zero, there is a large discrepancy in size between the sum set and difference set of $G$.  In principle, this should lead to a bound on $\alpha$, especially in view of standard inequalities relating sum sets and difference sets, such as
$$ |A-B| \leq { |A+B|^3 \over |A| |B| }.$$
(A summary of such inequalities can be found in \cite{ruzsa}).
However, these arguments (which are mostly graph-theoretical) do not seem to adapt well to the Kakeya application, because we only are working with a subset $G$ of $A_0 \times A_1$ rather than all of $A_0 \times A_1$.

To overcome this problem, Bourgain adapted a recent argument of T. Gowers which allows one to pass from arithmetic information on a subset of a Cartesian product to arithmetic information on a full Cartesian product.  A typical result is:

{\bf Theorem.}  Let $A, B$ be finite subsets of a torsion free abelian group with cardinality at most $N$, and suppose that there exists a set $G \subset A \times B$ of cardinality at least $\alpha N^2$ such that the sum set $\{ a+b: (a,b) \in G\}$ has cardinality at most $N$.  Then there exists subsets $A'$, $B'$ of $A$, $B$ respectively such that $A'-B'$ has cardinality at most $\alpha^{-13} N$ and $A'$, $B'$ have cardinality at least $\alpha^9 N$.

Roughly speaking, this theorem states that if most of $A+B$ is contained in a small set, then by refining $A$ and $B$ slightly, one can make {\it all} of $A-B$ contained in a small set also. Such results are reminiscent of standard combinatorial theorems concerning the size of sum and difference sets, but the innovation in Gowers' arguments is that the control on $A'$ and $B'$ is polynomial in $\alpha$.  (Previous combinatorial techniques gave bounds which were exponential or worse, which is not sufficient for Kakeya applications).

Recently \cite{kt:kakeya} N. Katz and the author have obtained the bound $(4n+3)/7$, by using control on both the sum-set $A_0+A_1$ and the variant $A_0+2A_1$, which corresponds to the slice $A_{2/3}$.  

These results have remarkably elementary proofs.  Apart from some randomization arguments, the proofs rely mainly on standard combinatorial tools such as the pigeonhole principle and Cauchy-Schwarz inequality, as well as basic arithmetic facts such as
$$ a+b = c+d \iff a-d = c-b,$$
$$ a-b = (a-b') - (a'-b') + (a'-b), \hbox{ and}$$
$$ a_0 + 2b_0 = a_1 + 2b_1, b'_0 = b'_1 \implies a_1 - b'_1 = 2(a_0 + b_0) - 2b_1 - (a_0 + b'_0).$$
Further progress has been made by pursuing these methods, though it seems that we are still quite far from a full resolution of the Kakeya problem, and some new ideas are almost certainly needed.

One possibility may be that one would have to use combinatorial estimates on {\it product sets} in addition to sum sets and difference sets, since one has control of $\{ a+tb: (a,b) \in G\}$ for all $t \in [0,1]$.  Discrete versions of such estimates exist; for instance, G. Elekes has recently shown the bound
$$ \max(|A \cdot A|, |A + A|) \geq C^{-1} |A|^{5/4}$$
for all finite sets of integers $A$.  However, these bounds do not adapt well to the continuous Kakeya setting because of the difficulty in discretizing both addition and multiplication simultaneously.  A good test problem in this setting is the {\it Erd\"os ring problem} - determine whether there exists a (Borel) sub-ring of $\R$ with Hausdorff dimension exactly 1/2.  This problem is known to be connected with the $\beta$-set problem and the Falconer distance set problem.

Interestingly, the Kakeya problem is also connected to another aspect of arithmetic combinatorics, namely that of locating arithmetic progressions in sparse sets.  (A famous instance of this is an old conjecture of Erd\"os that the primes contain infinitely many arithmetic progressions of arbitrary length).  This difficulty arises in the Hausdorff dimension formulation of the Kakeya problem, and also in some more quantitative variants, because of the difficulty in selecting a ``good'' set of three slices in arithmetic progression in which to run the above argument.  The combinatorial tools developed for that problem by W. Gowers and others may well have further applications to the Kakeya problem in the future.

\section{Applications to the Fourier transform}

Historically, the first applications of the Kakeya problem to analysis arose in the study of Fourier summation in the 1970s.

If $f$ is a test function on $\R^n$, we can define the Fourier transform $\hat f$ by 
$$ \hat f(\xi) := \int_{\R^n} e^{-2\pi i x \cdot \xi} f(x)\ dx.$$
One then has the inversion formula
$$ f(x) = \int_{\R^n} e^{2\pi i x \cdot \xi} \hat f(\xi)\ d\xi.$$
Now suppose that $f$ is a more general function, such as a function in the Lebesgue space $L^p(\R^n)$.  The Fourier inversion formula still holds true in the sense of distributions, but one is interested in more quantitative convergence statements.  Specifically, we could ask whether the partial Fourier integrals
$$ S_R f(x) := \int_{|\xi| \leq R} e^{2\pi i x \cdot \xi} \hat f(\xi)\ d\xi$$
converge to $f$ in (say) $L^p$ norm (the pointwise convergence question is also interesting, say for $L^2$ functions $f$, but this seems extremely difficult to show in two and higher dimensions.  In one dimension this was proven in a famous paper by L. Carleson).  By the uniform boundedness principle, this is equivalent to asking whether the linear operators $S_R$ are bounded on $L^p(\R^n)$ uniformly in $R$.  By scale invariance it suffices to show this for $S_1$:
$$ \| S_1 f\|_{L^p(\R^n)} \leq C \|f\|_{L^p(\R^n)}.$$

The operator $S_1$ is known as the \emph{disk multiplier}, because of the formula
$$ \widehat{S_1 f} = \chi_B \hat f$$
where $B$ is the unit disk in $\R^n$.  In one dimension it is a classical result of Riesz that this operator is bounded on every $L^p$, $1 < p < \infty$, and so Fourier integrals converge in $L^p$ norm.  (Indeed, in one dimension $S_1$ is essentially the Hilbert transform).  In higher dimensions, $S_1$ is bounded in $L^2$ thanks to Plancherel's theorem, however the behaviour in $L^p$ is more subtle.  One has an explicit kernel representation which roughly looks like
$$ S_1 f(x) \approx \int {e^{\pm i|x-y|} \over (1 + |x-y|)^{(n+1)/2}} f(y)\ dy;$$
to be more precise, one must use Bessel functions instead of $e^{\pm i|x-y|}$.  The kernel is only in $L^p$ when $p > {2n \over n+1}$, so it was natural by duality arguments to conjecture that $S_1$ was bounded when ${2n \over n+1} < p < {2n \over n-1}$.  In 1971, however, C. Fefferman proved the surprising

\proclaim{Theorem.}  If $n > 1$, then $S_1$ is unbounded on $L^p$ for every $p \neq 2$.  

In particular, one does not have $L^p$ convergence for the Fourier inversion formula in higher dimensions unless $p=2$.

Roughly speaking, the idea is as follows.  By duality it suffices to consider the case $p>2$.  Let $R$ be a large number, and let $T$ be a cylindrical tube in $\R^n$ with length $R$ and radius $\sqrt{R}$, and oriented in some direction $\omega_T$.  Let $\psi_T$ be a bump function adapted to  the tube $T$, and let $\tilde T$ be a shift of $T$ by $2R$ units in the $\omega_T$ direction.  Then a computation shows that
$$ |S_1 ( e^{2\pi i \omega_T \cdot x} \psi_T(x) )| \approx 1$$
for all $x \in \tilde T$.

\begin{figure}[htbp] \centering
\ \psfig{figure=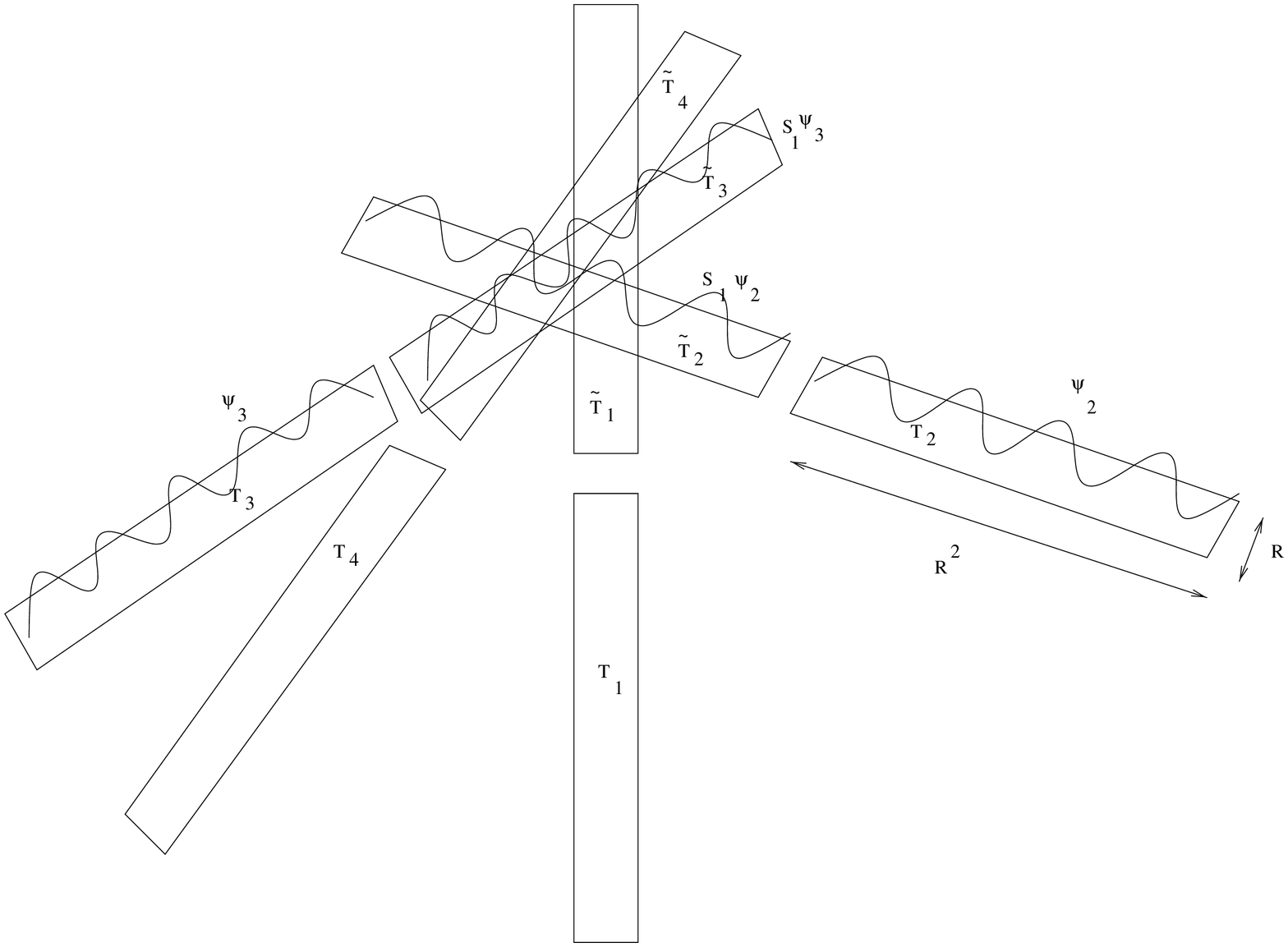,height=4in,width=6in}
\caption{Four tubes $T_i$, their shifts $\tilde T_i$, and the wave packets $\psi_i$.  The interference between the functions $S_1 \psi_i$ will cause the $L^p$
norm to be large when $p > 2$.
 }
\end{figure}

To exploit this computation, one uses the Besicovitch construction to find a collection $\{T\}$ of tubes as above, which are disjoint, but whose shifts $\tilde T$ have significant overlap.  More precisely, we assume that
$$ |\bigcup_T \tilde T| \leq \frac{1}{K} \sum_T |T|$$
for some $K$ which grows in $R$ (the standard construction gives $K \sim \log(R)/\log\log(R)$.  Then we consider the function
$$ f(x) = \sum_T \epsilon_T e^{2\pi i \omega_T \cdot x} \psi_T(x),$$
where $\epsilon_T = \pm 1$ are randomized signs.  Using Khinchin's inequality (which roughly states that one has the formula $|\sum_T \epsilon_T f_T| \sim (\sum_T |f_T|^2)^{1/2}$ with very high probability), one can eventually compute that
$$ \| S_1 f \|_p \geq C^{-1} K^{{1 \over 2}({1 \over 2} - {1 \over p})} \|f\|_p.$$
Since $K$ is unbounded, we thus see that $S_1$ is unbounded.

Fefferman's theorem is an example of how a geometric construction can be used to show the unboundedness of various oscillatory integral operators.  The point is while the action of these operators on general functions is rather complicated, their action on ``wave packets'' such as $e^{2\pi i \omega_T \cdot x} \psi_T(x)$ is fairly easy to analyze.  One can then generate a large class of functions to test the operator on by superimposing several wave packets together, and possibly randomizing the co-efficients to simplify the computation.

The counter-example provided by the Besicovitch construction is very weak (only growing logarithmically in the scale $R$), and can be eliminated if one mollifies the disk multiplier slightly.  For instance, the counterexample does not prohibit the slightly smoother {\it Bochner-Riesz operator} $S_1^\eps$, defined by
$$ \widehat{S_1^\eps(\xi) f} = (1 - |\xi|)^\eps \chi_B(\xi) \hat f(\xi)$$
from being bounded for $\eps > 0$, because the analogous computation gives
$$ |S_1^\eps ( e^{2\pi i \omega_T \cdot x} \psi_T(x) )| \approx R^{-\eps}$$
for all $x \in \tilde T$.  Indeed, {\it the Bochner-Riesz conjecture} asserts
that $S_1^\eps$ is indeed bounded on $L^p$ for all $\eps > 0$
and $2n/(n+1) < p < 2n/(n-1)$.  (For other values of $p$ one needs $\eps > n|{1 \over p} - {1 \over 2}| - {1 \over 2}$).  This conjecture was proven by L. Carleson and P. Sj\"olin in 1972 in two dimensions, but the higher-dimensional problem is quite challenging, and only partial progress has been made so far.  This conjecture would imply that the partial Fourier integrals will converge in $L^p$ if one uses a C\'esaro summation method (such as the Fejer summation method, which corresponds to $\eps = 1$).

The Bochner-Riesz conjecture would be disproved if one could find a collection of disjoint tubes $T$ for which the compression factor $K$ had some power dependence on $R$ as opposed to logarithmic, i.e. if $K \geq C^{-1} R^\eps$ for some $\eps > 0$.  A more precise statement is known, namely that the failure of the Kakeya conjecture would imply the failure of the Bochner-Riesz conjecture.  (More succinctly, Bochner-Riesz implies Kakeya).

In 1991, Bourgain introduced a method in which these types of implications could be reversed, so that progress on the Kakeya problem would (for instance) imply progress on the Bochner-Riesz conjecture.  The key observation is that every function can be decomposed into a linear combination of wave packets, by applying standard cutoffs both in physical space (by pointwise multiplication) and in frequency space (using the Fourier transform).  After applying the Bochner-Riesz operator to the wave packets individually, one has to re-assemble the wave packets and obtain estimates for the sum.  Kakeya estimates play an important role in this, since the wave packets are essentially supported on tubes; however this is not the full story since these packets also carry some oscillation, and one must develop tools to deal with the possible cancellation between wave packets.  The known techniques to deal with this cancellation - mostly based on $L^2$ methods - are imperfect, so that even if one had a complete solution to the Kakeya conjecture, one could not then completely solve the Bochner-Riesz conjecture.  Nevertheless, the best known results on Bochner-Riesz (e.g. in three dimensions the conjecture is known \cite{tv:cone} for $p > 26/7$ or $p < 26/19$) have been obtained by utilizing the best quantitative estimates of Kakeya type that are known to date.

These techniques apply to a wide range of oscillatory integrals.  A typical question, the {\it (adjoint) restriction problem}, concerns Fourier transforms of measures.  Let $d\sigma$ be surface measure on (say) the unit sphere $S^{n-1}$.  The Fourier transform $\widehat{d\sigma}$ of this measure can be computed explicitly using Bessel functions, and decays like $|x|^{-(n-1)/2}$ at infinity.  In particular, it is in the class $L^p(\R^n)$ for all $p > 2n/(n-1)$.  The \emph{restriction conjecture} asserts that the same statement holds if $\widehat{d\sigma}$ is replaced by $\widehat{f d\sigma}$ for any bounded function $f$ on the sphere.  This question originally arose from studying the restriction phenomenon (that a Fourier transform of a rough function can be meaningfully restricted to a curved surface such as a sphere, but not to a flat surface like a hyperplane); it is also related to the question of obtaining $L^p$ estimates on eigenfunctions of the Laplacian on the torus (although the eigenfunction problem is far more difficult due to number theoretic issues), as well as obtaining $L^p$ estimates on solutions to dispersive PDE, as we shall see below.

The restriction conjecture is logically implied by the Bochner-Riesz conjecture, and is slightly easier to deal with technically.  It has essentially the same amount of progress as Bochner-Riesz; for instance it is completely solved in two dimensions, and is true \cite{tv:cone} for $p > 26/7$ in three dimensions.  One uses the same techniques, namely wave packet decomposition of the initial function $f$, Kakeya information, and $L^2$ estimates to handle the cancellation, in order to obtain these results.

There are an endless set of permutations on these types of oscillatory integral problems; more general phases and weights, square function and maximal estimates, more exotic function spaces, bilinear and multilinear variants, etc.  
There are some additional rescaling arguments available in the bilinear case, as well as some $L^2$-based estimates, but apart from this there are few effective tools known outside of Bourgain's wave packet analysis to attack these types of problems.  

One variant of the above problems comes from replacing Euclidean space by a curved manifold.  There are some interesting three-dimensional examples of C. Sogge and W. Minicozzi which show that the Kakeya conjecture can fail on such manifolds, which then implies the corresponding failure of oscillatory integral conjectures such as the natural analogue of Bochner-Riesz.  This may shed some light on the robustness of Kakeya estimates and their applications in variable co-efficient situations.  Certainly the arithmetic and geometric techniques used to currently attack Kakeya problems do not adapt well to curved space.

\section{Applications to linear dispersive PDE}

We have seen in the previous section that Kakeya type problems have application to oscillatory integrals.  One important type of oscillatory integral arises from the solutions to linear PDE.  For instance, the solution to the free Schr\"odinger equation
$$ iu_t - \Delta u = 0; \quad u(0,x) = f(x)$$
for $x \in \R^n$, $t \in \R$ is given by the explicit formula
$$ u(t,x) = \frac{1}{(4\pi t)^{n/2}} \int_{\R^n} e^{-i|x-y|^2/4t} f(y)\ dy$$
for all $t \neq 0$.  We shall restrict our attention to the free (homogeneous) equation for sake of exposition, but much of this discussion can be extended to the forced (inhomogeneous) problem
$$ iu_t - \Delta u = F(x,t); \quad u(0,x) = 0,$$
which is also of importance in applications.

One way to see the connection between the free Schr\"odinger equation and the problems in the previous section is to take space-time Fourier transforms
$$ \hat u(\tau,\xi) := \int_{\R^{n+1}} e^{-2\pi i (x \cdot \xi + t \tau)} u(t,x)\ dt dx.$$
The Schr\"odinger equation then becomes
$$ (-2\pi \tau + 4 \pi^2 |\xi|^2) \hat u(\tau,\xi) = 0; \quad \int_\R \hat u(\tau,\xi)\ d\tau = \hat f(\xi).$$
From these formulae one easily obtains that $\hat u$ must be the distribution
$$ \hat u(\tau,\xi) = \hat f(\xi) \delta(\tau - 2\pi |\xi|^2)$$
with $\delta$ being the Dirac delta.  In particular, $u$ is the inverse Fourier transform of a measure on the paraboloid $\tau = 2 \pi |\xi|^2$.

Thus the problem of estimating the size of solutions to the Schr\"odinger equation is closely related to the restriction problem mentioned in the previous section.  By using the work of P. Tomas and E. Stein on the restriction problem, R. Strichartz showed a number of estimates for these types of equations.  A typical one is
$$ \| u \|_{L^4_{t,x}(\R \times \R^2)} = (\int_{\R^{2+1}} |u(t,x)|^4\ dt dx)^{1/4} \leq C \| f\|_{L^2(\R^2)}$$
for all solutions to the two-dimensional Schr\"odinger equation with finite $L^2$ norm.  Estimates of this type are known as {\it Strichartz estimates}; they encapsulate certain smoothing properties of the solution (an $L^4$ function is smoother than an $L^2$ function) as well as decay properties in time (the spatial $L^4$ norm of $u$ must eventually go to zero in order to make the space-time $L^4$ norm finite).  They are essential in establishing the well-posedness and scattering theory of certain non-linear dispersive equations with low regularity initial data (see the next section), especially when combined with more elementary estimates such as $L^2$ conservation
$$ \| u(t) \|_{L^2(\R^n)} = \| f \|_{L^2(\R^n)} \hbox{ for all } t \in \R.$$

The above Strichartz estimate can be proven in a number of ways, either by a direct calculation using the Fourier transform, or by taking adjoints and using restriction theory, or by breaking $u$ and $f$ into wave packets and using some rudimentary Kakeya information.  It has many generalizations; in fact, any linear equation which has a conserved energy and whose solutions decay in time will enjoy some family of Strichartz estimates.

Strichartz estimates are one way of controlling the behaviour of linear PDE, but they are not the only such estimates.  For instance, consider the problem of whether the explicit formula for the solution $u(t,x)$ to the free Schr\"odinger equation actually converges back to the initial data $f(x)$ when $t \to 0$.  For test functions $f$ this is easily verified, but the issue is more delicate when $f$ is rough (cf. the discussion on the Fourier inversion formula above).  A typical assumption is that $f$ is only in a Sobolev space $H^s$ for some $s \in \R$, or in other words that
$$ \| f \|_{H^s} := \| (1 + \sqrt{-\Delta})^s f \|_{L^2}$$
is finite.  (Roughly speaking, a function is in $H^s$ if the first $s$ derivatives of the function are in $L^2$.  One can make sense of the operator $(1 + \sqrt{-\Delta})^s$ by using the Fourier transform).

One can ask the question of whether $u(t,x)$ converges pointwise a.e. to $f(x)$ as $t \to 0$.  This can be shown to be equivalent to estimating the maximal function
$$ \sup_{t > 0} |u(t,x)|$$
in some appropriate norm in terms of the $H^s$ norm of the initial data $f$.  It is fairly easy to show that a.e. convergence fails for $s < 1/4$ and holds for $s > 1/2$, but the intermediate region requires a deeper understanding of the Schr\"odinger equation, and is not completely settled in two and higher dimensions.  For instance, in two dimensions the best known positive result \cite{tv:cone} is for $s > 1/2 - 1/32$.  This argument relies on a bilinear estimate related to the restriction problem, which in turn is proven using wave packet decompositions and Kakeya estimates.  (One may view this type of argument as decomposing the Schr\"odinger wave into ``particles'', and using Kakeya estimates to control how often these particles collide with each other as time evolves).

\begin{figure}[htbp] \centering
\ \psfig{figure=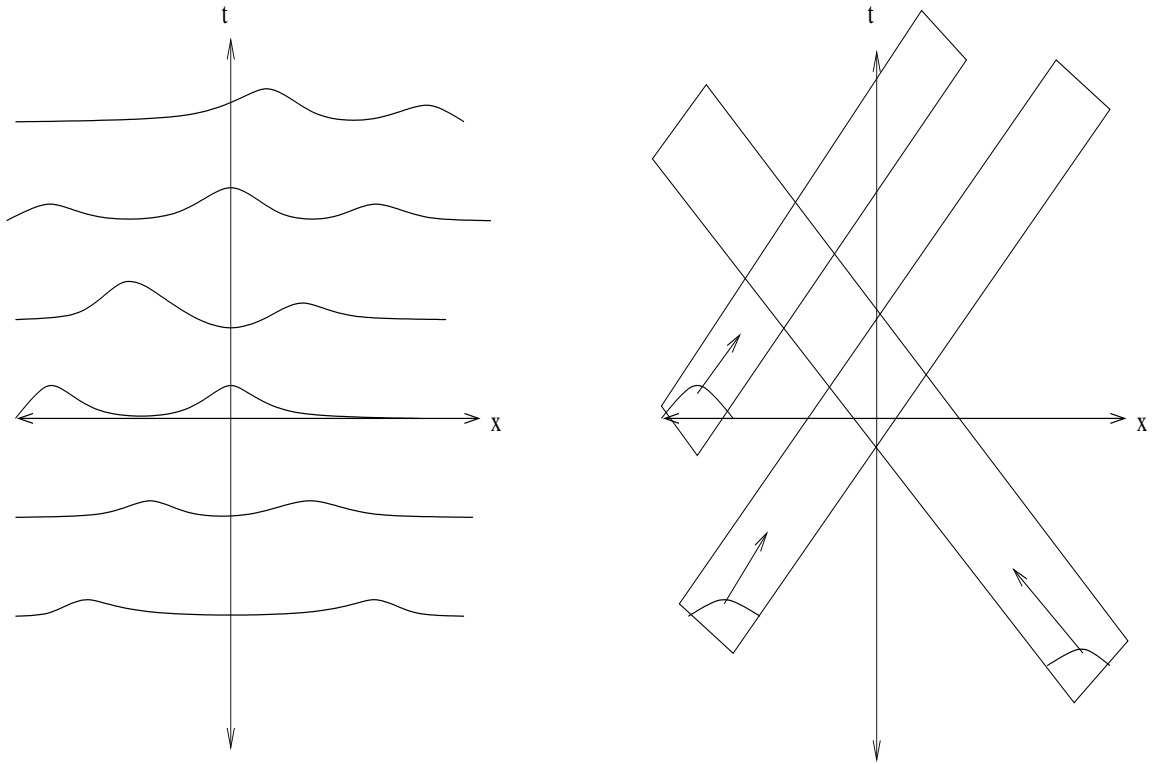,height=4in,width=6in}
\caption{A schematic depiction of how a wave (such as the one drawn on the
left) can be written as a superposition of
wave packets or ``photons''.  These objects are localized in space and
also have a localized direction, and different wave packets are essentially orthogonal.  There is no canonical way to perform this decomposition,
but one usually uses a combination of spatial cutoffs and cutoffs in 
Fourier space.
 }
\end{figure}

The quantitative behaviour of the free Schr\"odinger equation is only partially understood at present.  The Strichartz estimates are completely classified, but these estimates do not give complete control on the solution.  One basic issue (important in the study of energy concentration) is to understand the relationship between dispersion of $u(x,t)$, and the distribution of the Fourier transform $\hat u(x,t)$ on the paraboloid $\{ \tau = 2\pi |\xi|^2\}$.  Heuristically, one expects $u$ to disperse unless $\hat u(x,t)$ is concentrated in a small ``cap'' on the paraboloid, but obtaining a precise quantification of this phenomenon remains elusive, although estimates in this direction have been obtained and some interesting conjectures have been made.  It seems likely that more precise Kakeya estimates need to be developed before one has a satisfactory understanding of the size and dispersion of Schr\"odinger solutions.

There is a close parallel between the theory of the Schr\"odinger equation and that of the free wave equation
$$ \Box u(t,x) = 0; \quad u(0,x) = f(x), u_t(0,x) = 0$$
where $\Box = -\partial_t^2 + \Delta$.  To simplify the exposition we restrict to the case of zero initial velocity.

One can solve for $u$ explicitly using the formula
$$ u(t) = \cos(t \sqrt{-\Delta}) f$$
but this does not reveal much information about the size and distribution of $u$.  The Fourier transform of $u$ is supported on the light cone $\{ |\tau| = |\xi| \}$, which plays the same role as the paraboloid in the Schr\"odinger equation, although the cone has less curvature (so $u$ has less dispersion), and the cone has bounded slope (so $u$ has finite speed of propagation).

A general class of problem is the following: given size and regularity conditions on the initial data $f$ what type of size and regularity control does one obtain on the solution?

From integration by parts (or from the above explicit formula) one has energy conservation
$$ \int \frac{1}{2} |\nabla u(t,x)|^2 + \frac{1}{2} |u_t(t,x)|^2\ dx
= \int \frac{1}{2} |\nabla f(x)|^2$$
for all time $t$.  This conservation law, and its generalizations, show that $u$ has as much regularity as $f$ when measured $L^2$ based spaces.  However, $L^2$ control by itself does not reveal whether $u$ focuses or disperses.  To obtain better quantitative control on $u$ one needs other estimates, such as $L^p$ estimates.

The energy estimate is a \emph{fixed time} estimate; it controls the solution at a specified time $t$.  In the $L^p$ setting, fixed time estimates are very unfavorable, and require a lot of regularity for the initial data.  A typical estimate is the \emph{decay estimate}
$$ \| u(t) \|_{L^\infty(\R^n)} \leq C (1 + |t|)^{-(n-1)/2} \sum_{0 \leq k \leq s} \| \nabla^k f \|_{L^1(\R^n)}$$
whenever $s > (n+1)/2$ is an integer.  The necessity of this many derivatives is demonstrated by the \emph{focussing example}, in which the initial data is spread out near a sphere of radius $1$, and the solution $u$ focuses (with an extremely high $L^\infty$ norm) at the origin at time $t=1$.

However, one can obtain much better estimates, requiring far fewer derivatives, if one is willing to average locally in time.  The intuitive explanation for this is that it is difficult for a wave to maintain a focus point (which would generate a large $L^p$ norm for $p > 2$) for any length of time.  This phenomenon is known as \emph{local smoothing}.  Strichartz estimates are global smoothing estimates and thus fall into this category.  A typical such estimate is
$$ \| u \|_{L^4_{x,t}(\R^{3+1})} \leq C \| (\sqrt{-\Delta})^{1/2} f \|_{L^2_x(\R^3)}$$
in three spatial dimensions; this is the analogue of the Strichartz estimate for the two-dimensional Schr\"odinger equation mentioned earlier.  Without the averaging in time, one would require $3/4$ of a derivative on the right-hand side rather than $1/2$; this can be seen from Sobolev embedding.  These Strichartz estimates are usually proven by combining the energy and decay estimates with some orthogonality arguments.

However, even Strichartz estimates lose some regularity.  One may ask if there are $L^p$ estimates other than the energy estimate which do not lose any derivatives at all.  Unfortunately, even if one localizes in time and assumes $L^\infty$ control on the initial data, one still cannot do any better than $L^2$ control, as the following result of Wolff shows:

\proclaim{Theorem.}  If $n > 1$ and $p > 2$, then the estimate
$$ \| u \|_{L^p([1,2] \times \R^n)} \leq C \| f \|_{L^\infty(B(0,1))}$$
cannot hold for all bounded $f$ on the unit ball.

The argument proceeds similarly to Fefferman's disk multiplier argument. Let $\{ T \}$ be a collection of disjoint tubes arranged using the Besicovitch set construction as in Fefferman's argument, except that we rescale the tubes to have dimensions $1 \times R^{-1/2}$ rather than $R \times \sqrt{R}$.  On each of these tubes $T$ we place a ``wave train'', which is basically $e^{i R x \cdot \omega_T}$ times a bump function adapted to $T$.  Let $f$ be the sum of all these wave trains (although we may randomize the signs of these trains to simplify computations).

\begin{figure}[htbp] \centering
\ \psfig{figure=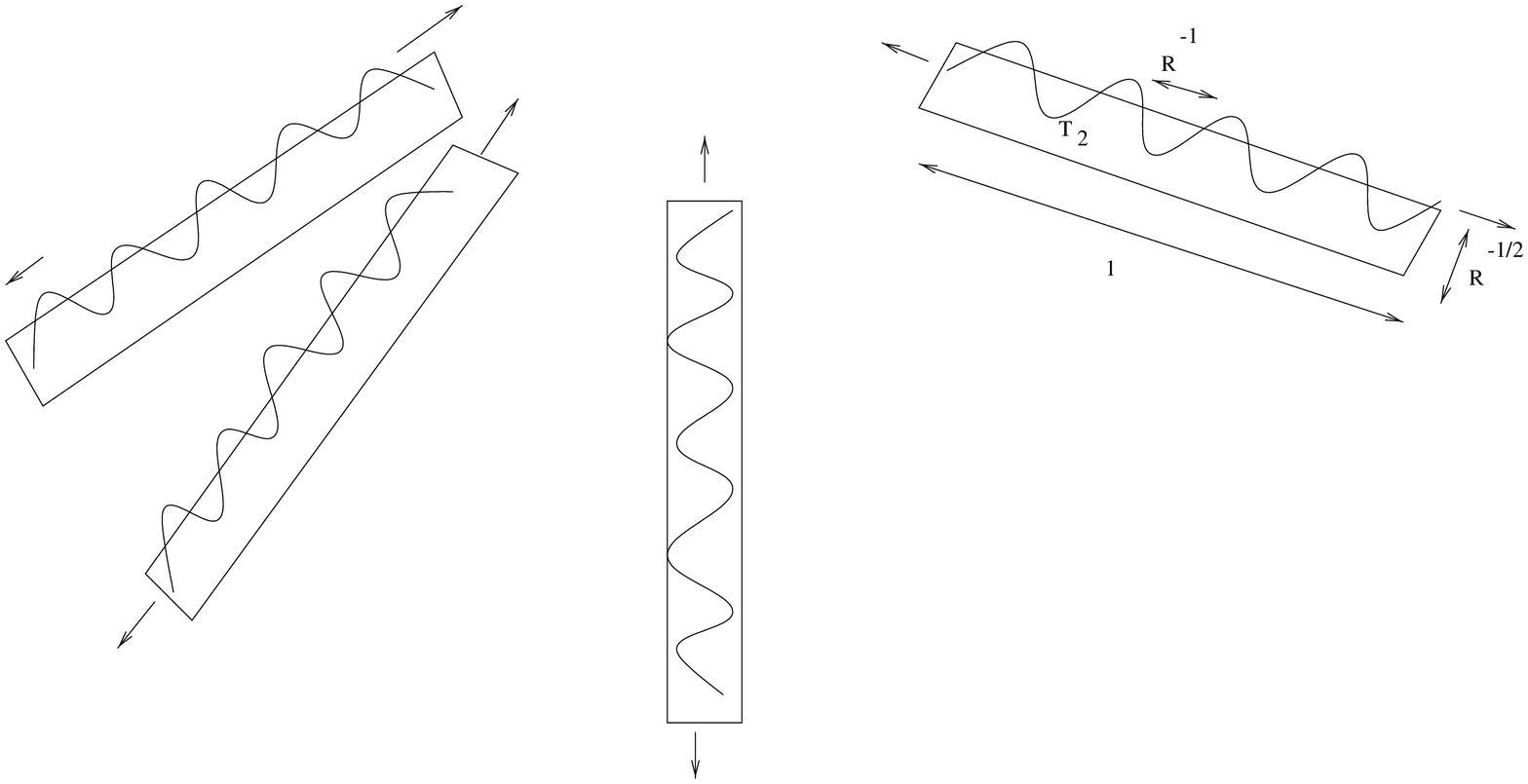,height=4in,width=6in}
\caption{A schematic depiction of the Wolff example at time zero.  Because of
the zero initial velocity, the wave trains will move in two opposite
directions.}
\end{figure}

\begin{figure}[htbp] \centering
\ \psfig{figure=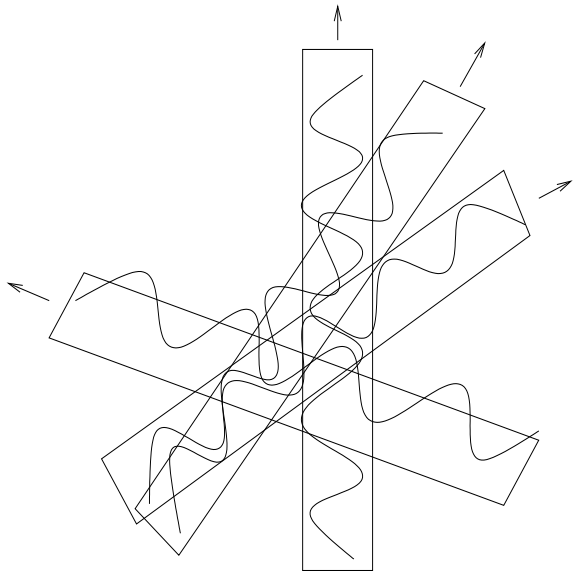,height=3in,width=3in}
\caption{The Wolff example at a time $1 \leq t \leq 2$; only the incoming wavetrains are shown.}
\end{figure}

At time zero, the function $f$ has low $L^\infty$ norm.  However, as time evolves, each wave train $T$ splits as the superposition of two pulses, one moving in the direction $\omega_T$, and the other in the direction $-\omega_T$.
For times $1 \leq t \leq 2$, a large portion of the wave train at $T$ now lives in the shifted tube $\tilde T$.  Because of the large overlap of these tubes, the $L^p$ norm of $u$ is large for all $1 \leq t \leq 2$; as with Fefferman's argument, it is about $K^{{1 \over 2}({1 \over 2} - {1 \over p})}$.  By letting $R \to \infty$ one can make $K$ unbounded, and this gives the Theorem.

Because the Besicovitch construction has a logarithmic compression rate, one could get around this obstruction by requiring an epsilon of regularity on the initial data.  The \emph{local smoothing conjecture} of C. Sogge asserts that no
further loss of regularity occurs, or more precisely that
$$ \| u \|_{L^p([1,2] \times \R^n)} \leq C_{p,\eps} \| (1 + \sqrt{-\Delta})^\eps f \|_{L^p(\R^n)}$$
for all $\eps > 0, n(\frac{1}{2}-\frac{1}{p}) - \frac{1}{2}$ and $2 \leq p \leq \infty$.  This conjecture is easy when $p=2$ or $p =\infty$; the most interesting case is when $p = 2n/(n-1)$.

The local smoothing conjecture is extremely strong, and would imply many of the known estimates on the wave equation.  It implies the Kakeya conjecture, for a counter-example to the Kakeya conjecture could be used to strengthen Wolff's argument to disprove the local smoothing conjecture.  This conjecture also implies the Bochner-Riesz conjecture; the idea is to write the Bochner-Riesz multiplier $S_1^\eps$ in terms of wave operators $\cos(t\sqrt{-\Delta})$ via a one-dimensional Fourier transform.  However, the conjecture is far from settled; even in two dimensions, the conjecture is completely proven for only $p > 74$ (due to T. Wolff), and at the critical exponent $p=4$ the conjecture is only known for $\eps > 1/8 - 1/88$ \cite{tv:cone}, \cite{Wo:cone}.  

There are several other wave equation estimates which are related to those discussed here.  An active area of research is to obtain good {\it bilinear} estimates on solutions to the wave equation, as opposed to the linear estimates described here; this will be discussed further in the next section.  

Since $u$ can be written in terms of circular averages of $f$, there is a also close relationship between wave equation estimates and estimates for circular means.  (Such circular means estimates can then be used, for instance, to make progress on the Falconer distance problem mentioned earlier).  There is also an extremely strong square function estimate conjectured for the wave equation which, if true, would imply the local smoothing, Bochner-Riesz, restriction and Kakeya conjectures.  It would also give estimates for other seemingly unrelated objects such as the helix convolution operator $f \mapsto f * d\sigma$, where $d\sigma$ is arclength measure on the helix $\{ (\cos t, \sin t, t): 0 \leq t \leq 2\pi\}$ in $\R^3$.  (The connection arises because the Fourier transform of $d\sigma$ is concentrated near the light cone).  These estimates are quite difficult, and the partial progress which has been made on them has proceeded via Kakeya estimates.  Although these deep wave equation estimates have not yet found significant applications, I am confident that they will do so in the near future.

\section{Applications to non-linear dispersive PDE}

Some of the estimates on linear PDE discussed in the previous section have been proven to be very useful in the analysis of non-linear perturbations of these PDE.  We shall restrict our discussion mostly to non-linear wave equations, although there has been parallel lines of research on similar equations such as non-linear Schr\"odinger and KdV.  A detailed discussion of these topics can be found in the recent survey \cite{ks}.

The idea that estimates on linear PDE lead to control of non-linear perturbations is an old one, dating back at least to Picard.  To illustrate the basic idea, consider the problem of solving the non-linear equation
$$ 5 x + \eps \sin(x) = b,$$
for $x \in \R$, where $\eps$, $b$ are given, and $\eps$ is small.  Such an equation cannot be solved explicitly, however one can show that a solution exists, is unique, and depends analytically on the parameters $\eps$ and $b$ provided that $\eps$ is sufficiently small.  To see this, simply rewrite the equation as
$$ x = 5^{-1} b + \eps 5^{-1} \sin(x).$$
Since $5^{-1} \sin(x)$ is Lipschitz, the right-hand side is a contraction if $\eps$ is sufficiently small, and the claim then follows from the contraction mapping principle.  More concretely, one constructs the solution $x$ as the limit of the iteration
$$ x^{(n+1)} := 5^{-1} b + \eps 5^{-1} \sin(x^{(n)})$$
with $x^{(0)} = 0$ (say).  An essentially equivalent approach would be to expand $x$ as a power series in $\eps$.

This perturbative approach extends to non-linear evolution equations provided that the non-linear equation is sufficiently well approximated by some linear or otherwise easily solvable equation ($5x=b$, in the above example).   This method is also very useful for demonstrating that qualitative phenomena associated to a special equation are robust under perturbations of that equation; for instance, one can use these techniques show that many non-linear wave equations obey the same type of properties (decay, resolution into plane waves, finite speed of propagation, energy conservation, few interactions between high and low frequencies, etc.) that the free wave equation enjoys.

The iteration method gives existence, uniqueness, and analytic dependence on the data; these properties are collectively referred to as {\it analytic well-posedness}.  In order to make the iteration converge, one often needs to localize in time, so that one only obtains local analytic well-posedness. 

Iterative techniques do not extend well to strongly non-linear equations such as Navier-Stokes, but have been quite successful when applied to dispersive equations such as the Korteweg-de Vries (KdV) equation and its variants, as well as non-linear Schr\"odinger and wave equations.  The latter category includes such equations as the wave maps and Yang-Mills, as well as the more difficult quasi-linear equations such as the Einstein equations.

A typical example is the semi-linear wave equation
$$ -u_{tt}(t,x) + \Delta u(t,x) = -u^3(t,x); u(0,x) = \eps f(x); u_t(0,x) = 0$$
in three spatial dimensions, where we set the initial velocity to zero for sake of exposition, and the initial position to be small.  This is one of the simplest non-linear perturbations of the free wave equation.  

The ODE $u_{tt} = u^3$ blows up in finite time, as one can see from such solutions as $u(t) = \sqrt{2}/(t-t_0)$.  So one does not always expect a global solution, especially for large initial data.  However, the Laplacian term in $\Box$ creates a dispersion effect, which can counteract this ODE blowup mechanism.  Heuristically, if $\eps$ is sufficiently small, then the dispersive effect of the Laplacian should overcome the blowup effect of the non-linearity, and one should obtain global existence.

We can rewrite this equation in integral form as
$$ u(t) = \eps \cos(t \sqrt{-\Delta}) f + \Box^{-1} u^3(t)$$
where $\Box^{-1} F$ is the unique solution $u$ to the inhomogeneous wave equation
$$ \Box u = F; u(0,x) = u_t(0,x) = 0.$$
An explicit formula for $\Box^{-1}$ is given by Duhamel's formula:
$$ \Box^{-1} F(t) = \int_0^t \frac{\sin((t-s)\sqrt{-\Delta})}{\sqrt{-\Delta}}
F(s)\ ds.$$

Consider the first few iterates of this equation.  The first iterate would be the free solution $u_0 := \eps \cos(t \sqrt{-\Delta}) f$.  The second iterate would be obtained by inserting the free solution into the non-linearity, to obtain a better approximation $u_1 := u_0 + \Box^{-1} u_0^3$.  And so forth.  It is clear that in order to make these approximations converge we must obtain estimates on such expressions as $\Box^{-1} u_0^3$.  The $L^p$ norms are suited for this task, and by using the Strichartz estimates for the free wave equation (together with variants for the inhomogeneous wave equation), 
one can indeed show analytic global well-posedness 
for data $f$ in the Sobolev space $\dot H^{1/2}$ for sufficiently small $\eps$, as the above integral equation turns out to be a contraction in a small ball in $L^4_{x,t}$.  In fact, one can refine the analysis and also obtain scattering results (so that the non-linear solution eventually converges to a linear solution).  The space $\dot H^{1/2}$ is not chosen arbitrarily; it is the unique Sobolev space $H^s$ which is invariant under the natural scaling associated with the above equation, and is referred to as the \emph{critical regularity} for this equation.  Regularity smoother 
(resp. rougher) than critical is referred to as \emph{subcritical} (resp. 
\emph{supercritical}) respectively.  As a general rule, local well-posedness 
is relatively straightforward to establish for sub-critical regularities, quite delicate (and often unknown) for critical regularities, and for super-critical regularities the techniques described above break down completely.  For super-critical regularities one does not expect well-posedness, let alone analytic well-posedness, although it is quite plausible that solutions will exist and be unique for ``generic'' choices of supercritical data.

There are several motivations for obtaining well-posedness results near the critical regularity, beyond the general benefit of testing, refining, and understanding the techniques used to prove these results.  Firstly, the time of existence given by such results tends to grow as one approaches the critical regularity, and in the critical case one often gets global existence automatically for small data.  Secondly, if blowup does indeed occur, the near-critical well-posedness theory often gives very precise control on the nature of the blowup (e.g. the energy must concentrate at a point).  In certain cases, such control 
can then be combined with other facts (such as conservation laws) to ensure that blowup cannot actually occur.  Finally, one can often refine the well-posedness analysis to obtain good control on various aspects of the solution, such as the transfer of energy from low frequency modes to high frequency modes or vice versa; this can also be used to prevent blowup and demonstrate global well-posedness even in cases where the energy is infinite.

For semi-linear equations such as the one described above, a satisfactory theory has been constructed using the Strichartz estimates as the main technical tool.  The situation is slightly less satisfactory with equations which are still semi-linear, but have derivatives in the non-linearity.  A typical example is the wave map equation for the sphere.  Wave maps are maps $\phi: \R^{n+1} \to S^m$ to a sphere which obey the non-linear wave equation
$$ \Box \phi = -\phi ( |\nabla \phi|^2 - |\phi_t|^2 ).$$
An important unsolved problem is whether wave maps with smooth data and finite energy remain smooth for all time in the two-dimensional case $n=2$.  (The result is false for $n > 2$ and true for $n=1$, but the $n=2$ case is far more delicate because the energy norm is critical in this case).

If one tries to apply the above techniques, one is quickly faced with the need to estimate bilinear expressions such as
$$ Q(\phi,\psi) := \nabla \phi \cdot \nabla \psi - \phi_t \cdot \psi_t.$$
Most bilinear combinations of waves are strongest when $\phi$, $\psi$ are both parallel travelling waves, e.g. $\phi = c_1 e^{ik_1(x_1 - t)}$,
$\psi = c_2 e^{ik_2(x_1-t)}$.  However, the above expression actually vanishes in this case.  Because of this, the above expression is known as a \emph{null form}.

One could estimate this null form by the triangle and H\"older inequalities, then using linear estimates such as the Strichartz estimates.  However, this method does not take advantage of the cancellation in the null form and requires a fair amount of regularity on the solution (about half a derivative more than the critical regularity).  In recent years it has been realized that in order to estimate these null forms accurately one needs genuinely bilinear estimates that take advantage of the null cancellation.

As was observed already, null forms vanish when $\phi$ and $\psi$ have the parallel frequencies.  Thus, one expects the largest contribution to a null form to come from interactions between waves which are transverse (i.e. their frequency vectors make an angle comparable to 1).  The best tool we currently have to control such interactions are bilinear $L^2$ estimates.  A typical such estimate is
$$ \| Q(\phi, \psi) \|_{L^2_{x,t}(\R^{3+1})} \leq C \| f \|_{H^1(\R^3)} \|g\|_{H^2(\R^3)}$$
where $\phi$, $\psi$ are free waves with initial position $f$, $g$ respectively and initial velocity zero, and $Q$ is a null form.  Strichartz methods can give this estimate, but with the right-hand side replaced by the symmetrized variant $\|f\|_{H^{3/2}} \|g\|_{H^{3/2}}$, which is less useful for applications.  This particular estimate, together with some variants, was used by S. Klainerman and M. Machedon to demonstrate global well-posedness of Yang-Mills and Maxwell-Klein-Gordon fields for finite energy data in three spatial dimensions.

These $L^2$ estimates can be proven by taking Fourier transforms and using Plancherel's theorem.  However, one can also prove these estimates using Kakeya methods.  Roughly speaking, the idea is as follows.  We shall induct on a scale parameter $R$; more precisely, we shall show that for each $R > 0$ we have the bound
$$ \| Q(\phi, \psi) \|_{L^2_{x,t}(Q_R)} \leq C \| f \|_{H^1(\R^3)} \|g\|_{H^2(\R^3)}$$
for all space-time cubes of side-length $R$.  Letting $R \to \infty$ we obtain the result.

Now fix $R$.  We decompose the free wave $\phi$ into wave packets, each of which has a somewhat localized position and velocity.  (This can be achieved by using cutoff functions in physical space to localize the initial data $\phi(0,x)$, then using Fourier cutoffs to localize the direction of travel and the frequency magnitude).  This decomposes $\phi$ into a superposition of functions $\sum_T \phi_T$ which are essentially supported on tubes near light rays (see Figure 9).  (The size of these tubes depends on $R$ and of the frequency magnitude of $\phi$.  For instance, if $\phi$ only contains frequencies of magnitude comparable to 1, then these tubes will have width about $\sqrt{R}$).

\begin{figure}[htbp] \centering
\ \psfig{figure=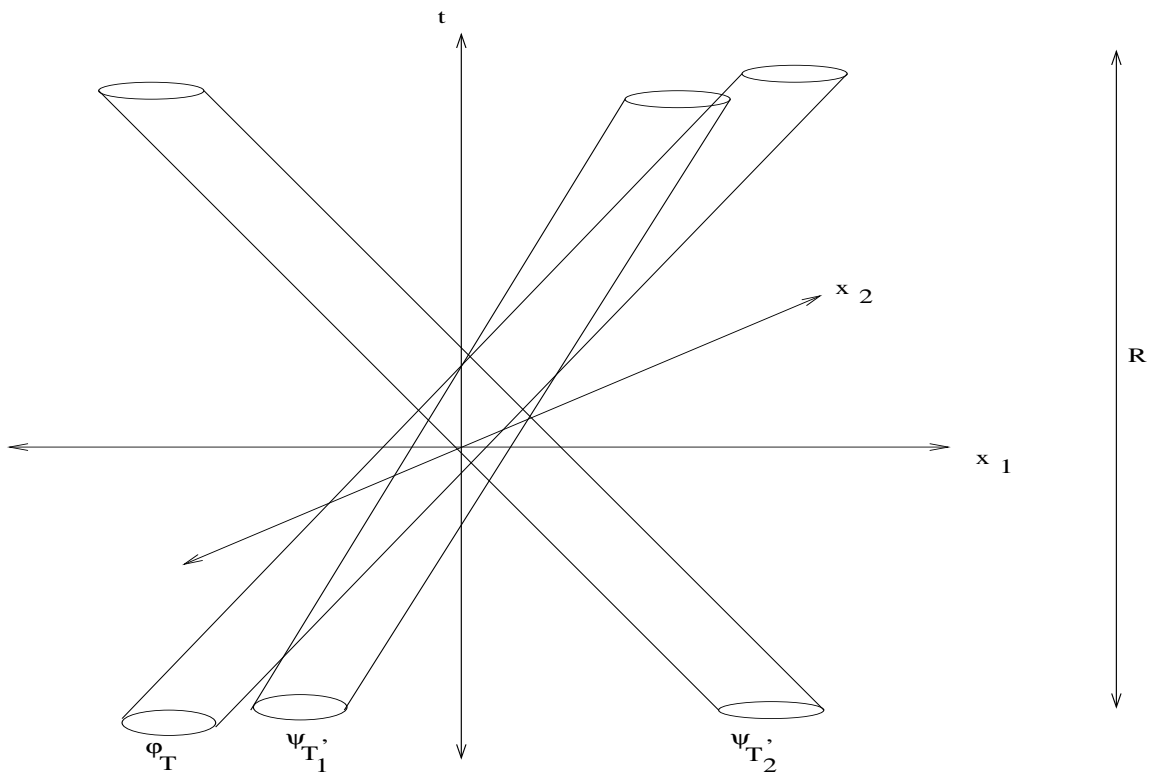,height=4in,width=6in}
\caption{The wave packet decomposition.  In this case $\phi$ consists
of a single wave packet $\phi_T$, while $\psi$ is the superposition of two
wave packets $\psi_{T'_1}$ and $\psi_{T'_2}$.  The parallel interaction between
$\phi_T$ and $\phi_{T'_1}$ would be small because of the cancellation in
the null form.}
\end{figure}
 
We can decompose $\psi = \sum_{T'} \psi_{T'}$ similarly.  We thus can split the original null form $Q(\phi, \psi)$ into a large number of smaller null form interactions $Q(\phi_T, \psi_{T'})$.

One can divide these interactions into {\it parallel interactions}, in which $T$ and $T'$ are parallel or nearly parallel, and {\it transverse interactions}, in which $T$ and $T'$ intersect with angle comparable to 1.  (If $T$ and $T'$ do not intersect at all, then $Q(\phi_T, \psi_{T'})$ is negligible).  The exact dividing line between what is parallel and what is transverse will depend on $R$ and on the magnitude of the frequencies of $\phi$ and $\psi$.  The parallel interactions are small because the null form vanishes when applied to parallel plane waves, so we concentrate on the transverse case.  Now we use the basic geometric fact that two transverse tubes of length $R$ and thickness $\sqrt{R}$ can only intersect in a cube of size roughly $\sqrt{R}$.   To exploit this, we partition the original cube $Q_R$ into smaller cubes $q$ of side-length $\sqrt{R}$, and observe that any pair $(T, T')$ of tubes can only interact within a finite number of such smaller cubes.  If one applies the induction hypothesis to each of $q$ individually, and then sums up over $q$, one can (essentially) close the induction and recover the $L^2$ estimate at scale $R$ from the $L^2$ estimates at scale $\sqrt{R}$.  (One has to use some orthogonality properties of the $\phi_T$ and $\psi_{T'}$ to properly do this).

\begin{figure}[htbp] \centering
\ \psfig{figure=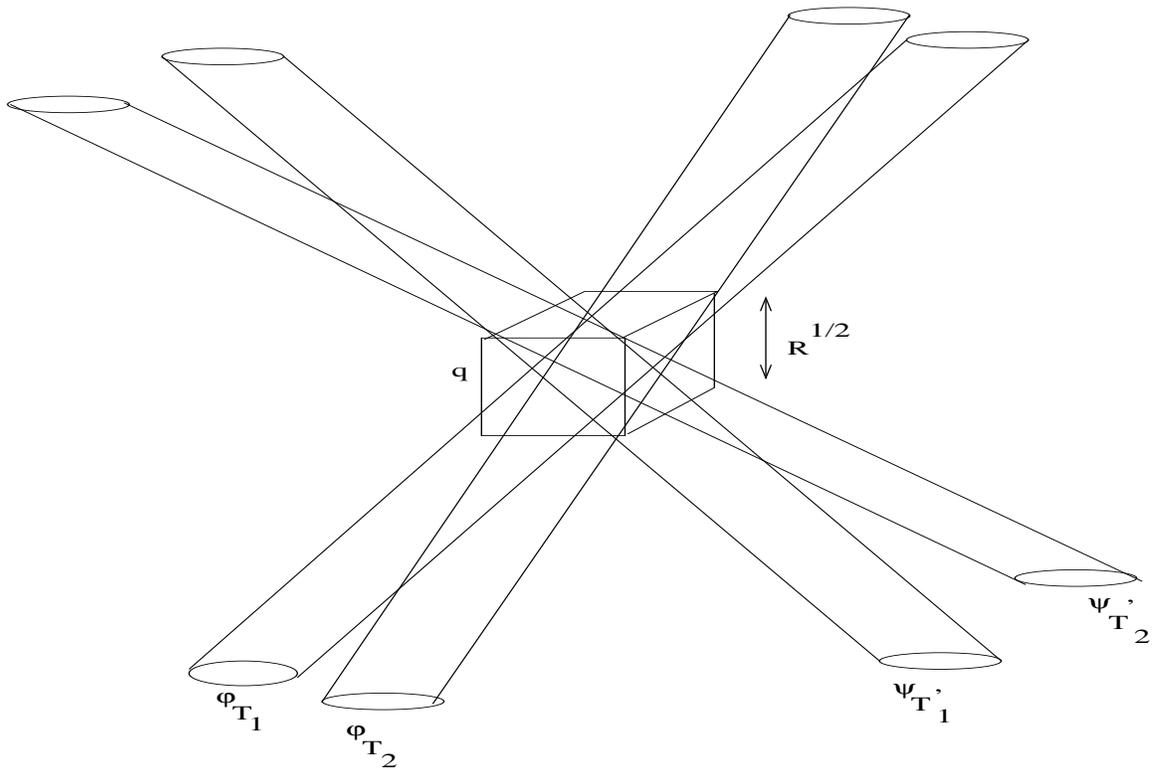,height=4in,width=6in}
\caption{The transverse interactions occuring at a single cube $q$ of
side-length $\sqrt{R}$.  By applying the induction hypothesis at this scale,
and using orthogonality of wave packets to sum over all cubes $q$, one can
close the induction.}
\end{figure}

This induction on scales idea is mainly due to T. Wolff, who also extended the above argument to cover a near-optimal range of $L^p$ estimates in addition to $L^2$ estimates.  These estimates and techniques are very recent, and have not yet found direct application to non-linear equations, but I am confident that they will do so in the near future.

There are several active areas of research in non-linear wave equations.  One important problem is to understand what happens at critical regularities.  The iteration techniques described above seem to just barely fail to control the solution in many of these cases, and some new ideas are needed.  Another problem is to extend the local existence theory given by iteration techniques to global existence.  If there is a finite conserved quantity which controls the local time of existence, then this is straightforward, but a current area of research (by J. Bourgain, M. Keel, J. Colliander, G. Staffilani, H. Takaoka, the author, and others) is to see what can be done if the conserved quantity is infinite.  Another area, in which much progress has recently been made, is to extend the semi-linear theory to {\it quasi-linear} equations, in which the $\Box$ operator is replaced by a variable co-efficient operator depending on the solution $\phi$.  A typical equation is
$$ (\partial_t^2 - \sum_{i,j} g_{ij}(\phi) \partial_i \partial_j) \phi = Q(\phi,\nabla^{-1} \phi),$$
where $g_{ij}$ is a smooth perturbation of the Kronecker delta $\delta_{ij}$, and $\nabla^{-1}$ is some pseudo-differential operator of order -1.  This is a quasi-linear equation of type similar to the Yang-Mills or Maxwell-Klein-Gordon equations.  So far, progress on this problem (by H. Bahouri, J. Y. Chemin, D. Tataru and S. Klainerman) has come from adapting Strichartz estimates to variable co-efficient situations when the metric is quite rough.  A current area of research is to also adapt the bilinear estimates to the rough co-efficient setting.  For this the Fourier transform is no longer as efficient a tool, but it seems that physical space techniques such as Kakeya methods are more robust and should be able to handle this situation.

\section{Acknowledgements}

I thank Nets Katz for some helpful references. The author is a Clay long-term Prize Fellow, and is supported by grants from the Sloan and Packard foundations.  This article is based on a talk given at the Clay Millennium event.

\end{document}